\documentclass[english,polish,12pt]{amsart}
\usepackage{babel}
\usepackage[T1]{fontenc}
\usepackage[cp1250]{inputenc}
\usepackage[centertags]{amsmath}
\usepackage{amsfonts}
\usepackage{amsthm}
\usepackage{newlfont}
\usepackage{amscd}
\usepackage{amsgen}
\usepackage{pb-diagram}
\usepackage{lamsarrow}
\usepackage{pb-lams}
\newlength{\defbaselineskip} 
\theoremstyle{plain}
\newtheorem{thm}{Theorem}[section] \newtheorem{cor}{Corollary}[section] \newtheorem{lem}[section]{Lemma} 
\newtheorem{lemm}{Lemma}[section]
\newtheorem{theo}{Theorem}
\theoremstyle{remark} \newtheorem{rem}{Remark}[section]
\theoremstyle{definition} \newtheorem{defi}{Definition}[section]
\theoremstyle{definition} \newtheorem{ex}{Example}[section] %
\newtheorem{constr}{Construction}[section]
 \numberwithin{equation}{section}
\begin{document}
\selectlanguage{english}
\title{Equations of log del Pezzo surfaces of index $\leq 2$ }
\author{Grzegorz Kapustka and Michał Kapustka}
\thanks{The work was done during the authors stay at the University of Liverpool
supported by the Marie Curie program.} \thanks{MSC-class: 14J26;
14J28, 14J10}
\begin{abstract}
Del Pezzo surfaces over $\mathbb{C}$ with log terminal
singularities of index $\leq 2$ were classified by Alekseev and
Nikulin. In this paper, for each of these surfaces, we find an
appropriate morphism to projective space. These morphisms enable
us to describe log del Pezzo surfaces of index $\leq 2$ by
equations in some weighted projective space. Obtained results give
a natural completion of similar results of Du Val, and Hidaka and
Watanabe, describing del Pezzo surfaces of index 1.
\end{abstract}
\maketitle
\section*{Introduction}
A log del Pezzo surface is a complex surface with quotient
singularities whose anti-canonical divisor is ample. The index of
a log del Pezzo surface $Z$ is the minimal number $i\in
\mathbb{N}$ such that the divisor $i K_Z$ is a Cartier divisor.

Log del Pezzo surfaces of index 1 have been studied by Du Val in
\cite{DVa,DVb}, Demazure in \cite{D}, and finally classified by
Hidaka and Watanabe in \cite{HW}. Hidaka and Watanabe gave two
classifications, the first one is in terms of blowing ups and
blowing downs, whereas the second one is in terms of equations in
some weighted projective space.

In presented paper, we study log del Pezzo surfaces of index $\leq
2$. These surfaces were classified by Nikulin and Alekseev in
terms of diagrams (graphs) of exceptional curves on an appropriate
resolution of singularities (see \cite{AN}). This is a counterpart
of the first classification of Hidaka and Watanabe. Our aim is to
extend their second description to this larger family of log del
Pezzo surfaces. That is, we are interested in finding simple
equations in some weighted projective space describing them.

The main results of the paper (Theorems \ref{weighted g2},
\ref{weighted g3}, \ref{g4} and \ref{higher genus}) may be
formulated as the following theorem.
\begin{theo}If $Z$ is a log del Pezzo surface of index $2$ then
we have one of the following cases:
\begin{enumerate}
 \item[a)] $K_Z^2 =1$ and $Z$ is isomorphic to a hypersurface of
 degree 6 in $\mathbb{P}(1,1,2,3)$
 \item[b)] $K_Z^2 =2$ and $Z$ is isomorphic either to a
 hypersurface of degree 8 in $\mathbb{P}(1,1,4,4)$ or to a
 hypersurface of degree 4 in $\mathbb{P}(1,1,1,2)$
 \item[c)] $K_Z^2 =3$ and $Z$ is isomorphic to a complete
 intersection in $\mathbb{P}(1,1,1,1,2)$ of two hypersurfaces of
 degree 2 and 3 respectively
 \item[d)] $K_Z^2 \geq 4$ and $Z$ can be embedded in
 $\mathbb{P}(1^{K_Z^2+1},2)$ in such a way that its ideal is
 generated by elements of degree 2 and 3. Moreover $Z$ is then:
 \begin{enumerate}
 \item[(i)] either isomorphic to a cone over a rational normal quartic
 curve
 \item[(ii)] or contained in a threefold cone isomorphic to a
 cone over a Veronese embedding of a rational normal scrollar
 surface or $\mathbb{P}^2$.
 \end{enumerate}
\end{enumerate}
\end{theo}

The main idea of the proof is to study embeddings of K3 surfaces
associated by Alekseev and Nikulin (\cite{AN}) to each considered
del Pezzo surface.

In section \ref{intro} we revise the construction of associated
$K3$ surfaces and collect some preliminary results concerning log
del Pezzo surfaces and their resolutions.

In section \ref{sec1} we introduce a linear system $|D_g|$ on the
constructed K3 surface and then, we use results of Saint-Donat
(see \cite{SD}) to describe the morphism given by it. The morphism
descent to a mapping from the log del Pezzo surface to the
quotient of the projective space by a hyperplane symmetry. In most
cases, our description is related to the morphism given by the
anti-canonical divisor.

We consider each of the cases from the theorem in a separate
section. In each of these sections beside proving the appropriate
part we make a deeper analysis of the obtained examples.

Section \ref{sec2} deals with surfaces $Z$ for which $g=K_Z ^2 +1
=2$. The aim of that section is to give a complete description of
log del Pezzo surfaces having above property. The first theorem
(Thm. \ref{M}) gives a description in terms of coverings. More
precisely we find a 2:1 morphism onto a quadric cone in
$\mathbb{P}^3$ (that is onto the weighted projective space
$\mathbb{P}(1,1,2)$) and describe its branch locus. The equations
follow then easily (Thm. \ref{weighted g2}).

In this section, obtained equations occur as natural deformations
(and some kind of completion) of analogous equations described in
\cite{HW} for Gorenstein log del Pezzo surfaces (i. e. of index
1).

For instance Du Val, in \cite{DVa, DVb}, introduces Du Val
singularities as arising from the double covering of a quadric
cone branched over the intersection of this cone with a smooth
cubic surface not passing through its vertex. This surface is
Gorenstein and log del Pezzo. We obtain the cases where the cubic
passes through the vertex as models of log del Pezzo surfaces of
index 2. This surfaces have a singularity over the vertex of the
cone that is not Du Val. It is a singularity of the so-called type
$K_n$.

The third section deals with the case when $K_Z ^2=2$. In this
section we are adapting to the case $g=K_Z ^2 +1=3$ the same
technics which had been used earlier for $g=K_Z ^2 +1=2$. The
difference is that the morphism $|D_g|$ is no more 2:1. The
description looses a little bit of its clarity because it not free
from checking singularities of resulting equations.

Here we also get a natural generalization of the analogous case
from \cite{HW}. But there is one example that does not arise as a
natural deformation of any of the surfaces described there.
Elements from this new family happens to be the only ones that has
more than one singularity of index 2 (they have exactly 2
singularities of type $K_1$).

The next section gathers all remaining cases. After proving the
main theorem we describe each of the possibilities from its
assertion in a more precise way in context of the diagrams of
exceptional curves from \cite{AN}. However, for $K_Z ^2 \geq 5$
the description is not so precise as earlier. Instead we use
constructions given in \cite{AN} to give a complete list of
families of log del Pezzo surfaces of index 2 with $g\geq 6$. We
can moreover find a good description for some of the examples from
this list, separately from the rest, by realizing them as quintic
hypersurfaces in $\mathbb{P}(1,1,1,4)$.

The last section shows an application of these results from a
different point of view. With the help of \cite{AN}, we get the
list of all configurations of singularities that can occur in each
of the surfaces described throughout all the paper. For example,
we can list all possible configurations of ADE singularities on
symmetric quartics in $\mathbb{P}^3$, extending the description of
\cite{HW} with the cases where the quartic passes through the
center of symmetry. By studying intersections of a quadric cone
with a cubic in $\mathbb{P}^3$, we can do the same for symmetric
sextics in $\mathbb{P}^2$. \setcounter{section}{-1}

\section*{Acknowledgements}
We would like to express our gratitude to Prof. V. Nikulin for his
enormous help and fruitful discussions during our work on this
paper.

\section{Notation and basic tools}\label{intro} In this section, we recall the notation
and basic facts about log del Pezzo surfaces of index $\leq 2$,
contained in \cite{AN}, that we will use in the paper.

First observe that the only singularities of index $\leq 2$ on log
del Pezzo surfaces are Du Val (of index 1) or of the type $K_n$.
\begin{defi}A singularity on a surface is called of type
$K_n$ if the exceptional divisor of its minimal resolution is one
of the following configurations of $n$ curves:
 \item[--] a rational curve of self intersection $-4$ for $n=1$
 \item[--] two rational curves of self intersection $-3$ cutting in one
 point for $n=2$
 \item[--] two rational curves of self intersection $-3$ connected
 by a chain of $-2$ curves for $n\geq 3$
\end{defi}
In other terms the diagram of exceptional curves in the minimal
resolution is the following:
\begin{itemize}
\item[(i)] \setlength{\unitlength}{0.5mm}
 \begin{picture}(10,10)(-2.5,-1)
 \put(0,0){\circle{3}}
 \put(0,0){\circle{1.5}}
 \end{picture} for n=1
\item[(ii)] \setlength{\unitlength}{0.5mm}
 \begin{picture}(20,10)(-2.5,-1)
 \put(0,0){\circle{3}} \put(0.3,0){\makebox(0,0){\tiny{x}}}
 \put(10,0){\circle{3}}\put(10.3,0){\makebox(0,0){\tiny{x}}}
 \put(1.5,0){\line(1,0){7}}
 \end{picture}
 for n=2
\item[(iii)]\setlength{\unitlength}{0.5mm}
 $\underbrace{
 \begin{picture}(65,10)(7.5,9)
 \put(10.0,10.0){\circle{3}}\put(10.3,10){\makebox(0,0){\tiny{x}}}
 \put(20.0,10.0){\circle*{3}}
 \put(30.0,10.0){\circle*{3}}
 \put(60.0,10.0){\circle*{3}}
 \put(70.0,10.0){\circle{3}}\put(70.3,10){\makebox(0,0){\tiny{x}}}
 \put(11.5,10){\line(1,0){7}}
 \put(21.5,10){\line(1,0){7.0}}
 \put(31.5,10){\line(1,0){7.0}}
 \put(61.5,10){\line(1,0){7.0}}
 \put(51.5,10){\line(1,0){7}}
 \put(46,10){\makebox(0,0){\dots}}
\end{picture}}_n$ \hskip0.7cm for $n\geq 3$
\end{itemize}
Where a black vertex denotes a rational $-2$ curve a double
transparent vertex a rational $-4$ curve and a crossed transparent
vertex a rational $-3$ curve.

To work with log del Pezzo surfaces, we need to resolve their
singularities. In our context, the interesting resolution will not
be the minimal one, but the right resolution introduced in
\cite{AN}.
\begin{defi}
If $Z$ is a log del Pezzo surface of index $\leq 2$ its resolution
of singularity is called right if and only if:
 \item[(i)] it is a minimal resolution of the Du Val singularities
 \item[(ii)] the singularities of type $K_n$ are resolved by
 composing the minimal resolution with the blowing up of all
 points of intersection of curves in the exceptional divisor.
 The exceptional divisor takes the form:
\setlength{\unitlength}{0.5mm}
 $$\underbrace{
 \begin{picture}(65,10)(7.5,9)
 \put(10.0,10.0){\circle{3}}
 \put(20.0,10.0){\circle{3}}
 \put(30.0,10.0){\circle{3}}
 \put(60.0,10.0){\circle{3}}
 \put(70.0,10.0){\circle{3}}
 \put(10.0,10.0){\circle{1.5}}
 \put(30.0,10.0){\circle{1.5}}
 \put(70.0,10.0){\circle{1.5}}
 \put(11.5,10){\line(1,0){7}}
 \put(21.5,10){\line(1,0){7.0}}
 \put(31.5,10){\line(1,0){7.0}}
 \put(61.5,10){\line(1,0){7.0}}
 \put(51.5,10){\line(1,0){7}}
 \put(46,10){\makebox(0,0){\dots}}
\end{picture}}_{2n-1}$$
where transparent vertices denote $-1$ curves.
\end{defi}

\begin{defi}
A right DPN surface of elliptic type surface (for us shortly a DPN
surface) is a smooth surface that we can obtain as a right
resolution of a log del Pezzo surface of index $\leq 2$.
\end{defi}
For basic properties and characterization of DPN surfaces we refer
the reader to \cite{AN}, where they were defined in a different
but equivalent way.

From now on, $Z$ will denote a log del Pezzo surface of index
$\leq 2$, on which we shall put additional conditions depending on
the subfamily we are going to work with.

We will denote by $g$ the number $K_Z ^2 +1$. This invariant plays
an important role in the classification of log del Pezzo surfaces
of index $\leq 2$ given in \cite{AN}. It has a geometrical
interpretation, as we can see below in the following construction.
\begin{constr}[see {\cite[sec 2.1]{AN}}]
Let $Z$ be a log del Pezzo surface of index $\leq 2$. Let $r:Y
\longrightarrow Z$ be its right resolution. $Y$ is then a right
DPN surface of elliptic type. Denote by $E_i$ (for $i \in
\{1,\dots, k\}$) the $(-4)$ exceptional curves of $r$. From
\cite[thm 1.5]{AN}, we know that in the system $|-2K_Y|$, there
exists a smooth curve $C$ of the form $C= C_g +E_1+E_2+ \dots
+E_k$ (in particular $C_g .E_i =0$), where $C_g \in |r ^{\ast }
(-2K_Z)| $ is a smooth curve of genus $g$. We can now take the
double covering $\pi : X \longrightarrow Y$ of $Y$ branched in the
curve $C$. By adjunction formula, $X$ is a K3 surface.
\end{constr}
\begin{rem}
In the above construction, we have a large choice of curves $C_g$,
each of them giving by this construction a different K3 surface
$X$. We will choose one of them and keep this notation as fixed
throughout the paper (however, in a suitable place, we will show
what changes if we make another choice).
\end{rem}

We can see from the construction that, any right DPN surface of
elliptic type can be obtained as a quotient of some K3 surface $X$
by an involution $\theta$ with fixed locus being a smooth curve of
the form $D_g +F_1 +\dots +F_k$, where $D_g$ is a curve of genus
$g\geq2$ and $F_i$ are rational. Such an involution on a K3
surface is called non-symplectic of elliptic type. The converse is
also true, namely if $X$ is a K3 surface and $\theta$ a
non-symplectic involution of elliptic type, then the quotient of
$X$ by $\theta$ is a right DPN surface of elliptic type.
\smallskip

An involution is called non-symplectic if its fixed locus is a
smooth curve. This type of involution is also characterized in the
following way. Let $\omega_X$ be any non-zero regular
2-dimensional differential form on $X$, then an involution
$\theta$ is non-symplectic if and only if $\theta
^{\ast}(\omega_X) = -\omega_X$. An involution on a K3 surface that
does not satisfy this condition has to be a symplectic involution,
that is $\theta ^{\ast}(\omega_X) = \omega_X$, and then its fixed
locus is a finite set of points.

\smallskip
Together with $g=K_Z ^2 +1$ in the classification contained in
\cite{AN} two more invariants are important. They are called $k$
and $\delta$ (see \cite[Sec. 2.3]{AN}).

\smallskip
The number $k$ is related to the type of singularities of index 2
appearing on the del Pezzo surface. It is the number of $-4$
curves on the corresponding DPN surface. It happens that in most
of the cases there is only one singularity of this index, it is by
definition a singularity $K_k$. The only exceptions are del Pezzo
surfaces with two $K_1$ singularities, they all form one family in
the description of \cite{AN} (in the Table 3 it has number 25) and
naturally have $k=2$ by definition.

\smallskip
The invariant $\delta$ takes only one of two values 0 or 1. The
$\delta=0$ if and only if the reduced curve $\pi ^{-1}(C_g + E_1 +
\dots + E_k)$ is divisible by 2 in the Picard lattice of the K3
surface $X$.

\medskip
Table 3 in \cite{AN} gives a complete list of possible diagrams of
exceptional curves on right DPN surfaces of elliptic type. It is
very long, so for lack of space we refer to the original paper.
However, we should recall some basic facts about how to understand
this list.

As earlier, black vertices denote $-2$ curves, transparent ones
$-1$ curves and double transparent ones $-4$ curves. Each of them
corresponds to a family of DPN surfaces that we construct by
blowing up $\mathbb{P}^2$, $\mathbb{P}^1 \times \mathbb{P}^1$ or
some $\mathbb{F}_n$ to obtain a configuration of exceptional
curves (curves with negative self-intersection) as in the diagram.

We denote the diagram of exceptional curves on a DPN surface $Y$
by $\Gamma(Y)$. The graph of black vertices in a diagram $\Gamma$
determines the type and the configuration of Du Val singularities
on the corresponding del Pezzo surfaces. It is called the Du Val
part of the diagram and is denoted by $\mathrm{Duv}(\Gamma)$. The
Table 3 in \cite{AN} does not contain all diagrams, each of them
generate a number of others corresponding to any subgraph of
$\mathrm{Duv}(\Gamma)$. In this way for any subgraph
$D\subset\mathrm{Duv}(\Gamma)$ there is a DPN surface with a
configuration of -2 curves described by $D$.

The double transparent vertices together with those transparent
ones that are connected with two double transparent ones form a
subgraph called the logarithmic part of the diagram $\Gamma$ and
denoted by $\mathrm{Log}(\Gamma)$.

To obtain a log del Pezzo surface from a DPN surface with a
diagram $\Gamma$ of exceptional curves, we need to contract the
curves from $\mathrm{Log}(\Gamma)$ and $\mathrm{Duv}(\Gamma)$.
These graphs describe the index 2 and the Du Val singularities of
the obtained del Pezzo surface respectively. We see that if
$\Gamma$ is a diagram from Table 3 in \cite{AN}, then for any
subgraph $D\subset\mathrm{Duv}(\Gamma)$ there is a del Pezzo
surface with a configuration of Du Val singularities described by
$D$.

We know also that the curves from a diagram of exceptional curves
of a DPN surface $Y$ represent all exceptional curves on $Y$.
Except the cases where $\mathrm{Pic}(Y)$ is generated by $\leq 2$
elements (this happens only when $Y$ is isomorphic to
$\mathbb{P}^2$, $\mathbb{P}^1 \times \mathbb{P}^1$ or some
$\mathbb{F}_n$), these exceptional curves generate
$\mathrm{Pic(Y)}=\mathrm{Num}(Y)$. This last equation follows from
the fact that $Y$ is obtained from $\mathbb{P}^2$, $\mathbb{P}^1
\times \mathbb{P}^1$ or some $\mathbb{F}_n$ by a sequence of
blowing ups. The diagram lifts up to a configuration of curves on
the K3 surface $X$ in such a way that the pre-images of double
transparent and transparent vertices are smooth rational curves on
$X$ and the pre-images of black vertices split up to two disjoint
components each of them being a smooth rational curve on $X$. The
obtained curves are all rational curves on $X$.

Observe (see \cite[Rem 3.2]{AN}) that it is easy to write down the
curve $C_g$ as numerically equivalent to a combination of vertices
of the diagram, as we know the intersection numbers with each of
these vertices. Indeed if $V$ is a curve represented by a vertex
of the diagram, then:
\begin{itemize}
\item $C_g . V = 0$ if $V$ is represented by a double transparent
vertex or a black vertex,
 \item $C_g . V =2-i$ if $V$ is represented by a transparent vertex
and $i$ is the number of intersection points of $V$ with curves
represented by double transparent vertices.
\end{itemize}

\begin{section}{General construction}\label{sec1}

Let us first do a general construction that we will use throughout
all the paper. We fix the notation from the discussion in section
\ref{intro}. Let $D_g \subset X$ be the reduction of the divisor
$\pi ^{\ast} C_g$. Observe that the curve $D_g$ is isomorphic to
$C_g$ and that $\pi ^{\ast} C_g = 2 D_g$

The aim of the construction is to find a natural and easy to
describe morphism from $Z$ to a weighted projective space. We will
use the morphism given by $|D_g |$ on the K3 surface $X$. First we
need to compute some invariants to construct a suitable diagram.
\begin{lemm} \label{dimensions of systems} The following formulas hold
\begin{itemize}\leftskip=-0.5cm
 \item[1)] $\dim |C_g|=3g-3$
 \item[2)] $C_g^2=4g-4$
 \item[3)] $\dim |D_g|=g$
 \item[4)] $D_{g}^2=2g-2$
\end{itemize}
\end{lemm}
\begin{proof}
Observe that $|2D_g|=|\pi ^{\ast}(C_g)|$ and use Riemann-Roch
theorem.
\end{proof}

On $X$, we have an involution determined by the covering $\pi$,
let us denote it by $i$. This $i$ fixes the curve $D_g$ so, acts
on its linear system. One obtains an involution on $H^0
(X,\mathcal{O}_X (D_g))$ and so on its projectivisation
$\mathbb{P}^{g}$. Such an involution, since it is not the
identity, is a symmetry $s$ with respect to two disjoint linear
subspaces generating $\mathbb{P}^{g}$ (because it is a symmetry in
the associated affine space). By definition, we have that:
\begin{equation}\label{symmetry of the image}
\varphi _{|D_g|}\circ i=s \circ \varphi _{|D_g|}.
\end{equation}
In particular the image $\varphi _{|D_g|}(X)$ is symmetric with
respect to $s$. Since the curve $D_g$ is contained in the fixed
part of $i$, its image $\varphi_{|D_g|}$ has to be contained in
one of the fixed subspaces. But this image spans a hyperplane $H$
in $\mathbb{P}^{g}$. Hence $s$ is a symmetry with respect to the
hyperplane $H$ and some point $O$. We can consider $\alpha :
\mathbb{P}^{g} \longrightarrow \mathbb{P}(1^g , 2)$ the quotient
morphism of the symmetry $s$. Now, we need only to find a morphism
to close the diagram as in the following lemma.
\begin{lemm}\label{Key} There exists a morphism $\varphi $ such
that the following diagram commutes:
\begin{equation}\label{diag}
\begin{diagram} \node{}
\node{X} \arrow{sw,l}{\pi} \arrow{e,t}{\varphi_{|D_g|}}
\node{\mathbb{P}^{g }} \arrow{s,l}{\alpha}  \\
\node{Y} \arrow{e,b}{r} \node{Z} \arrow{e,b}{\varphi}
\node{\mathbb{P}(1^g , 2)} \arrow{e,b}{\psi}
\node{\mathbb{P}^{\frac{g(g+1)}{2}}}
\end{diagram}
\end{equation}
Here $\psi$ denotes the natural embedding of $\mathbb{P}(1^g , 2)$
into $\mathbb{P}^{\frac{g(g+1)}{2}}$ given by all monomials of
weighted degree 2.
\end{lemm}

\begin{proof}
Observe that the curves $E_i$ are also invariant and are not in
the pre-image  $\varphi_{|D_g|}^*(H)$, so they have to map to the
point $O$. We know that the morphism $\alpha \circ
\varphi_{|D_g|}\colon X \longrightarrow \mathbb{P}(1^g , 2)$ is
invariant under $i$, so it factors through $\pi$ which is the
quotient morphism of the involution $i$. In a natural way, we have
obtained a morphism $\phi \colon Y \longrightarrow \mathbb{P}(1^g
, 2)$ making the diagram commutative. Observe that:
\begin{center}
$\pi ^* ((\psi \circ
\phi)^*(\mathcal{O}_{\mathbb{P}^{\frac{g(g+1)}{2}}}(1))) =(\psi
\circ \phi \circ \pi)^*(\mathcal{O}_{\mathbb{P}^{\frac{g(g+1)}{2}}}(1 )) =$\\
\medskip$ (\psi \circ \alpha \circ \varphi_{|D_g|})^*
(\mathcal{O}_{\mathbb{P}^{\frac{g(g+1)}{2}}}(1)) \subset |2
D_g|=|\pi ^* (C_g)|$.
\end{center}
By projection formula, we deduce that each element of the system
$(\psi \circ \phi)^*(\mathcal{O}_{\mathbb{P}^{\frac{g(g+1)}{2}}}
(1))$ is numerically, thus also linearly (see section
\ref{intro}), equivalent to $C_g$. This means that the morphism
$\psi \circ \phi$ is given by a subsystem of $|C_g|$. Hence $\psi
\circ \phi$ factors through the morphism $r$ that only contracts
curves $F$ such that $F.C_g=0$.
\end{proof}

\begin{rem}
Changing coordinates in $\mathbb{P}^{g}$ to $x_0 ,x_1,..,x_{g} $,
we can assume $H$ to be $x_{g}=0$ and $O$ the point $(0:0: ...
:0:1)$. From now on, we make this assumption.
\end{rem}
\begin{rem}
If in Lemma \ref{Key} we use another curve $C_g\in|-2K_Z|$, the
morphism $\varphi$ will stay unchanged. It is just the branch
locus of the quotient by the symmetry $s$ that changes. It simply
becomes another smooth hypersurface of weighted degree 2 in
$\mathbb{P}(1^g , 2)$.
\end{rem}
\end{section}
\begin{section}{Del Pezzo surfaces with $g=K_Z ^2 +1 =2$}\label{sec2}

We can now pass to the description of the first family of log del
Pezzo surfaces of index $\leq 2$. From now on, in this section,
$Z$ is a del Pezzo surface with $g=2$ (i. e. $K_{Z}^2= 1$). We can
formulate the main theorem.

\begin{thm}\label{M}\ \\
\textbf{\emph{1)}} Let $Z$ be a log del Pezzo surface with $g=K_Z
^2 +1 =2$. The map given by the linear system $|-2K_Z |$ is a
finite $ 2:1$ morphism onto a quadratic cone $Q$ in
$\mathbb{P}^3$. The branch divisor $R$ of this morphism is reduced
and is the intersection of the cone $Q$ with a smooth cubic (the
morphism is also branched in the vertex of the cone). Moreover $R$
is not the intersection of $Q$ with three planes meeting in a
point lying on this cone.\\
\textbf{\emph{2)}} Conversely, the normal double covering of a
quadratic cone $Q'\subset{P^3}$, bran\-ched over such an
intersection $R'$ (it is enough to assume that the cubic is smooth
in the vertex of the cone), is a log del Pezzo surface of index
$\leq 2$ with $g=2$.
\end{thm}
 Let us first prove the following lemma describing the
singularities that can occur on a symmetric plane sextic which
will appear in this section later.
\begin{lemm}\label{symmetric sextic}
Let $K$ be a reduced plane sextic, which is symmetric with respect
to a line $H$ and a point $O$, and not containing the symmetry
line $H$, and with at most a double singularity in the center of
symmetry $O$. Then $K$ can have only simple singularities or it
the sum of three tangent conics. Moreover, if $K$ passes through
$O$, it has at $O$ a singularity of type $A_n$ with $n$ being an
odd number.
\end{lemm}
\begin{proof}
Let us denote the symmetry with respect to $O$ and $H$ in
$\mathbb{P}^2$ by $s$. Choosing appropriate coordinates, we can
assume that \[s \colon \mathbb{P}^2 \ni(x:y:z)\mapsto (x:y:-z)\in
\mathbb{P}^2.\]

We first consider the case when $K$ passes through $O$. Suppose
that $K$ is singular in some point $P\not= O$. Then it is also
singular in its symmetric adjoint $P'=s(P)$ with the same
multiplicity. Take the line $L=PP'$. If $L$ is not a component of
$K$, then from the B\'ezout theorem the singularities in $P$ and
$P'$ have multiplicity at most 2. If $L$ is a component of $K$,
then $K$ is the sum of $L$ and some quintic $K_5$, such that
$L\not \subset K_5$. From B\'ezout theorem the singularities in
the points $P$ and $P'$ of $K_5$ have multiplicity at most 2 and
$L$ is not any of their tangent cones. This means, that the
tangent cones to $K$ at $P$ and $P'$ consist of at least two
lines. From the characterization of Du Val singularities (see
\cite{BPV}), we deduce that $K $ has in $P$ and $P'$ simple
singularities.
\medskip

As we have assumed the singularity of $K$ at the point $O$ is a
double singularity hence is of some type $A_n$. To determine this
$n$ consider the double covering $\pi_2\colon \widetilde{X}
\longrightarrow \mathbb{P}^2$ branched over the sextic $K$. The
surface $\widetilde{X}$ is then a K3 surface with Du Val
singularities, and a singularity $A_n$ in the inverse image of
$O$. We can see that $\widetilde{X}$ is given in $\mathbb{P}
(1,1,1,3)$ with coordinates $x,y,z,t$ by the equation $t^2
=F(x,y,z)$, where $F$ is the equation of $K$ in $\mathbb{P}^2$ .
Consider the involution $\mathbb{P} (1,1,1,3) \ni(x:y:z:t) \mapsto
(x:y:-z:t)\in \mathbb{P} (1,1,1,3)$, we can induce it to
$\widetilde{X}$. In this way, we obtain an involution
$\widetilde{i}$ on $\widetilde{X}$. Observe that $\pi_2 \circ
\widetilde{i} =s \circ \pi_2$. Consider the minimal resolution $X$
of $\widetilde{X}$. It is a K3 surface. The involution
$\widetilde{i}$ on $\widetilde{X}$ lifts up to an involution $i$
on $X$ in such a way that the inverse image on $X$ of the symmetry
point is a graph of curves of type $A_f$ which preserved by the
involution.

\smallskip
Let us look how $i$ acts on the graph $A_f$. This graph is
obtained from $\widetilde{X}$ as the graph of exceptional divisors
of successive blow ups and each blow up lifts the involution
$\widetilde{i}$. In this way all curves from the graph are
preserved by the involution $i$ and the first one is fixed. From
the description of involutions on K3 surfaces (see section
\ref{intro}) the fixed locus of $i$ is a smooth curve. In
particular $i$ cannot have isolated fixed points and no two fixed
curves from the graph intersect. Moreover as the graph $A_f$ of
exceptional curves does not intersect any fixed curve not from
this graph, every fixed point in this exceptional locus has to lie
on a fixed curve from the graph. Now it is enough to observe that,
if a curve is preserved by the involution, then this involution
induced on it is either the identity or has exactly two fixed
points. This means that every second curve from the chain $A_f$ is
fixed by $i$ and that the curves at both ends of the chain are
fixed. By this last argument, we see that we had to blow up an odd
number of times.

If $K$ does not pass through $O$, then we can use \cite{HW}.
\end{proof}

Let us pass to the proof of the theorem.
\begin{proof}[Proof of Theorem \ref{M}]\

\textbf{1)} We do the general construction described in the
previous section and keep the notation fixed there. We will first
describe the morphism given by $|C_g|= r ^{\ast}(|-2K_Z |)$ on the
DPN surface $Y$ (after that, we will factor it through $r$).

We have $g=2$ and $D_2$ is a smooth curve of genus 2 (hence it is
hyperelliptic) on a K3 surface. From \cite[5.1]{SD} this curve
gives a $2:1$ morphism to $\mathbb{P}^2$ branched over a reduced
sextic with Du Val singularities. We denote this curve by $R_6$.
This $R_6$ has to be symmetric with respect to $H$ and $O$,
because from the property \ref{symmetry of the image}, we have
$\varphi _{|D_2|}^{-1} (p)=i\circ\varphi _{|D_2|}^{-1} (s(p))$ and
$i$ is a bijection. Since the pullback to $X$ of the line
$H=\{z=0\}$ is $|D_2|$ (it is smooth of genus 2), the sextic $R_6$
has to cut $H$ transversely. In particular, it cannot have any
singularity on $H$. The diagram from Lemma \ref{Key} takes the
form.
$$
\begin{CD} X @>{\varphi_{|D_2|}}>2:1>\mathbb{P}^2 \\
@V{\pi }VV @ V {\widetilde{\alpha}} VV\\
 Y@ >>{\psi \circ \varphi \circ r } >\mathbb{P}^3\\
\end{CD} $$
Here $\widetilde{\alpha} := \psi \circ \alpha : \mathbb{P}^2
\longrightarrow \mathbb{P}^3$ is given by
$(x_0^2:x_1^2:x_0x_1:x_2^2)$, thus its image is a quadric cone
$Q$. Observe that the branch divisor of the map $\psi \circ
\varphi \circ r$ is the image by $\widetilde{\alpha}$ of the
sextic $R_6$.
\medskip

We claim that this branch divisor is a reduced intersection of the
cone $Q$ with a smooth cubic. \smallskip

Indeed, let $F$ be the equation of $R_6$ in $\mathbb{P}^2$. The
equation $F \circ s$ also defines $R_6$ so $F\circ s =F$ or
$F\circ s =-F$ since $F\circ s \circ s=F$. Now $F\circ s =-F$ is
impossible because the fixed line of the symmetry $s$ cannot be a
component of $R_6$ (as $R_6$ is reduced and symmetric of degree
6). We have $F\circ s =F$ so $F(x_0,x_1,x_2)=
G(x_0^2,x_1^2,x_0x_1, x_2^2)$ where $G(t,u,v,w)$ is a homogenous
equation of degree 3. We know that $R_6$ can have only a
singularity $A_n$ in the symmetry point $(0:0:1)$ so $G$ defines a
variety in $\mathbb{P}^3$ smooth in the point $(0:0:0:1)$.
Moreover, we can choose $G$ to define a smooth cubic in
$\mathbb{P}^3$, by adding a degree 3 polynomial divisible by
$(v^2-tu)$ (the generic element of this system will be good since
its intersection with the cone is reduced). \medskip

Coming back to the proof of the theorem, we have proved that the
morphism $\psi \circ \varphi \circ r$ is a morphism of degree 2
onto a cone $Q$ with branch divisor being the intersection of $Q$
with a smooth cubic $R$. Moreover from the proof of Lemma
\ref{Key} $( \psi \circ \varphi \circ r) ^{\ast } (
\mathcal{O}_{\mathbb{P}^3}(1) ) \subset |C_2|$. From Lemma
\ref{dimensions of systems} $\dim |C_2|=3$. This means that these
systems have the same dimension, hence are equal. Finally
$$\psi \circ \varphi \circ r =\varphi_{|C_g|}=\varphi_{r^{\ast}|-2K_Z|}.$$
This morphism contracts only the curves $F$ such that $F.C_g =0$.
Thus the morphism $\psi \circ \varphi$ is equal to
$\varphi_{|-2K_Z|}$, and it is a finite 2:1 morphism onto the cone
$Q$ branched in $Q\cap R$ and the vertex of the cone.

Since $\mathrm{Pic}(Q)$ has no torsion, the normal double covering
of $Q$ is uniquely determined by its branch divisor. Hence $Z$ is
isomorphic to the normal double covering of the cone $Q$ branched
in the intersection of $Q$ with a smooth cubic $R$. The last
assertion of 1) follows now from the fact that the normal double
covering of a cone branched over the intersection of this cone
with three hyperplanes meeting at a point on it, has an elliptic
singularity (see \cite{HW}). Thus it cannot be isomorphic to $Z$.
\bigskip

\textbf{2)} To prove the second part, we repeat the same reasoning
reconstructing it backwards. We will use analogous notation for
analogous objects but with added primes. Let $Z'$,$Q'$,$R'$ be
surfaces as in part 2) of the theorem. Let $H_1$ be a hyperplane
in $\mathbb{P}^3$ that cuts transversely the intersection curve of
the cone $Q'$ with the cubic $R'$. We can choose coordinates $t'$,
$u'$, $v'$, $w'$ in $P^3$ in such a way that $H_1=\{w'=0\}$ and
the point $(0:0:0:1)$ is the vertex of the cone $Q'$. Take the
normal double covering $\widetilde{\alpha}' \colon \mathbb{P}^2
\ni (x',y',z')\mapsto (x'^2, y'^2, y'z', z'^2) \in Q'\subset
\mathbb{P}^3$ of the quadric cone $Q'$ (that is $\mathbb{P}
(1,1,2)$ ) branched over the intersection of $Q'$ with $H_1$ (that
is over a degree 2 hypersurface in $\mathbb{P} (1,1,2)$). The
inverse image of this intersection is a reduced plane sextic $R_6
'$ symmetric with respect to the covering involution $s'$ (which
has to be the symmetry with respect to a line and a point). This
$R_6 '$ cuts transversely the fixed line $H'$ of this involution
and has a singularity of multiplicity $\leq 2$ in the center of
symmetry $O'$ (the inverse image of the vertex of the cone).

By Lemma \ref{symmetric sextic} $R_6 '$ must have Du Val
singularities. The desingularization of the double covering of
$\mathbb{P}^2$ branched over the sextic $R_6$ is a K3 surface
$X'$.
In the same way as in the proof of Lemma \ref{symmetric sextic},
we obtain an involution $i'$ on $X'$.

 The inverse image of the fixed line $H'$ is contained in the
fixed part of the involution $i'$ and since $R_6'$ cuts $H'$
transversely it is a curve of genus 2 (this is a double covering
of $\mathbb{P}^1$ branched in 6 points). Since the fixed locus of
an involution on a K3 surface has constant dimension, the fixed
locus of $i'$ is of codimension 1 in $X$. This implies that $i'$
has to be a non-symplectic involution of elliptic type. So using
\cite[thm 2.1]{AN} the quotient of $X'$ by this involution is a
right DPN surface of elliptic type. Contracting all exceptional
curves from $Log \Gamma (Y) $ and $Duv \Gamma (Y)$ (see \cite[thm
4.1]{AN}), one obtains a log del Pezzo surface $Z''$ of index 2
with $K_{Z''}^2=1$.

The normal surface $Z''$ is a double covering of the cone $Q'$
branched along the intersection of $Q'$ with the cubic $R'$, so it
is isomorphic to $Z'$. Thus $Z'$ is itself a log del Pezzo surface
of index $\leq2$ and $g=2$.
\end{proof}

\begin{rem} Du Val in \cite{DVa} considered the double covering
of a cone branched over the vertex and the intersection of this
cone with a smooth cubic not passing through the vertex. He
determines the types of singularities that arise in this way (Du
Val singularities). Here, we have a situation a bit more general.
In our context the cubic can pass through the vertex of the cone.
The case where it does pass gives rise to a little more
complicated singularity (of index 2) that we can determine exactly
in the following way.

As we know from \cite{AN} (see section \ref{intro}), the only
singularities of index 2 are of type $K_n$. Hence, we need to
determine $n$. We can do this as follows.

The minimal resolution of the cone $Q$ is the ruled surface
$\mathbb{F}_2$. We can easily compute that the pullback of the
cubic (passing trough the vertex of the cone) is then an element
$L\in |6f+2C|$ (where $C$ is the exceptional curve and $f$ a
fiber) that has neither multiple components nor $C$ as a
component. Now, if necessary by resolving the singularity of $L$
on $C$, we find $n$.
\end{rem}

The following theorem gives a description of log del Pezzo
surfaces of index $\leq 2$ and with $g=2$ in terms of equations in
weighted projective space.

\begin{thm}\label{weighted g2}\ \\
\textbf{\emph{1)}} Let $Z$ be a log del Pezzo surfaces of index
$\leq 2$ with $g=2$. Then $Z $ is isomorphic to a hypersurface in
$ \mathbb{P} (1,1,2,3) $ (with coordinates $x,y,z,t$) defined by
an equation F of degree 6 of the form
$$F(x,y,z,t)= t^2-G(x,y,z)=0$$
where the curve in $\mathbb{P}^2$  (with coordinates a,b,c) given
by the equation $G(a,b,c^2) $ is reduced, has at most a double
singularity in $(0,0,1)$ and is not the sum of three tangent
conics.\\
\textbf{\emph{2)}} Conversely each hypersurface $Z'$ in
$\mathbb{P} (1,1,2,3)$, defined by an equation as above is a log
del Pezzo surface of index $\leq 2$ with $g=2$.
\end{thm}
\begin{proof} First let us prove that the hypersurface in
$\mathbb{P} (1,1,2,3)$, defined by the equation $F$ is normal.

\smallskip
We know that $F$ defines a hypersurface with isolated
singularities. According to \cite{IF} it is enough check normality
in the singular points of $\mathbb{P} (1,1,2,3)$. But $\mathbb{P}
(1,1,2,3)$ has only two singularities $(0:0:0:1)$, $(0:0:1:0)$ and
since we have $F(0,0,0,1) \not= 0$, it is enough to check that
this hypersurface is normal in $(0:0:1:0)$. This point is obtained
as a quotient of the du Val singularity $t^2- G(x,y,1)=0$ (in the
affine piece $z=1$) by the action of $\mathbb{Z}_2$ (the symmetry
trough $(0,0,0)$). Since this action lifts to the universal
covering of the du Val singularity, the point $(0,0,1,0)$ is a
quotient singularity. Hence, according to \cite{B}, it is normal.
\medskip

To finish the proof, we use Theorem \ref{M} and the fact that the
normal double covering of a quadric cone is uniquely determined by
its branch divisor.
\end{proof}

\begin{rem}
 We can determine the main invariants ($g$, $k$, $\delta$) of the
obtained del Pezzo surface $Z'$ from the geometry of the sextic
($G=0$) in $\mathbb{P}^2$.

 We have $g=2$.

\smallskip If the sextic $\{G=0\}$ has a double point of
type $A_{f}$ in the symmetric point, then from Lemma
\ref{symmetric sextic} $f=2l+1$. Moreover, as we could see in the
proof of this lemma $l$ is the number of rational curves in the
fixed locus of the involution on the corresponding K3 surface $X$,
that is $l=k$.

\smallskip To compute the invariant $\delta $, we resolve the
singularities of the sextic $\{G=0\}$ and form a diagram of
exceptional curves on $X$. Next, we compute $\delta$ from the
definition. That is, we check how the obtained exceptional curves
intersect $D_2 + \pi^{-1}(E_1) + \dots + \pi^{-1} (E_k)$. Since we
have obtained the diagram of exceptional curves, we can also look
at \cite[table 3]{AN}.

\end{rem}

\end{section}


\begin{section}{Del Pezzo surfaces with $g=K_Z^2+1=3$}\label{sec3}
In this section $Z$ is a log del Pezzo surface of index $ \leq 2$
with $g=3$. As the main result, we get the following theorem which
gives us a description of all such surfaces.
\begin{thm}\label{g3}\ \\
\textbf{\emph{1)}} Let $Z$ be a log del Pezzo surface of index
$\leq 2$ and with $g=3$. Then $Z$ is isomorphic to one of the
following:
\begin{itemize}
\item[(a)] a hypersurface of degree 8 in $\mathbb{P}(1,1,4,4)$
(with coordinates $x$, $y$, $z$, $t$) given by an equation of the
form $t^2=G(x,y,z)$ where the equation $G(a,b,c^4)=0$ defines in
$\mathbb{P}^ 2$ (with coordinates $a$, $b$, $c$) a reduced octic
with
simple singularities 
which does not pass through the point $(0:0:1)$.

\item[ (b)] a hypersurface of degree 4 in $\mathbb{P} (1,1,1,2)$
given by an equation $F(x,y,z,t)=0$ where $ \{F(a,b,c,d^2) =0\}$
is a quartic in $\mathbb{P} ^4$ (with coordinates $a,b,c,d$) with
ADE singularities.
\end{itemize}
\textbf{\emph{2)}} Conversely, every surface described by
equations as in above cases is a log del Pezzo surface of index
$\leq 2$ with $g=3$.
\end{thm}

\proof Let $Z$ be a log del Pezzo surface of index $ \leq 2 $ with
$g=3$. Using the Lemma \ref{Key}, we get the commutative diagram
\ref{diag}. According to \cite[4.1]{SD}, we have two
possibilities. The morphism given by the linear system $|D_3|$ is:
\begin{itemize}
 \item either birational to a quartic with Du Val
singularities in $\mathbb{P} ^3$ (it is called the
non-hyperelliptic case)
 \item or 2:1 to a quadric in $\mathbb{P} ^3$
(it is called the hyperelliptic case)
\end{itemize}
Moreover, we know that the image $\varphi_{|D_3|}(X)$ is
symmetric. We shall separately consider the two cases.

\smallskip \textbf{The non-hyperelliptic case -} We can see in this
case, that the quartic does not have any singularity on the
hyperplane of symmetry. Then $\varphi$ from Lemma \ref{Key} is a
finite birational morphism to the quotient of the quartic
$\varphi_{|D_3|}(X)$ by the symmetry. Since this quartic is
irreducible, using a similar argument as in the proof of Theorem
\ref{M} we see that its equation is also symmetric. Hence, its
quotient is given (after a suitable change of coordinates) by an
equation in $\mathbb{P} (1,1,1,2)$ as in case  (b) of the theorem.
To prove that $\varphi$ is an isomorphism, we need only to remark
that the quotient of the quartic is normal. This is because the
only singularity that is not Du Val appears as a quotient of a Du
Val singularity, and this gives a quotient singularity (see proof
of Theorem \ref{weighted g2}).

\medskip
\textbf{The hyperelliptic case -} A priori, we have two sub-cases:

\smallskip
\textbf{The quadric is smooth}. Then the center of symmetry must
lie outside the quadric and taking the quotient, we get an
equation of degree 2 in $\mathbb{P} (1,1,1,2)$ not vanishing at
the point $(0:0:0:1)$, that is an equation of a surface isomorphic
to $\mathbb{P} ^2$. And $\varphi$ from the Lemma \ref{Key} is a
2:1 morphism to this surface. In this case $Z$ can have only Du
Val singularities (as the branch locus of $\varphi$ is locally
isomorphic to the branch locus of the morphism given by $|D_3|$,
i.e. has simple singularities). We can then obtain $Z$ as a double
covering of $\mathbb{P} ^2$ branched over a quartic with simple
singularities, that is case  (b) from the theorem, but, in fact we
will see in Remark \ref{impossible}, that this case never occurs.

\smallskip
\textbf{The quadric is a cone}. It is then isomorphic to
$\mathbb{P} (1,1,2)$. Let us denote this cone by $Q$. The image of
the curve $D_3$ is a hyperplane section of $Q$ hence a conic
$\widetilde{D}$ (it is smooth rational curve). But $D_3$ is a
smooth curve of genus 3, so the restriction of the map given by
$|D_3|$ to it, is a double covering of the conic $\widetilde{D}$
branched in 8 distinct points. Moreover the generic ray of the
cone $Q$ is the image of an elliptic curve, so the map, restricted
to it, is branched in 4 distinct points. This gives us the branch
locus to be a symmetric curve $R_8$ with only Du Val
singularities, which is the intersection of the cone with a
quartic (i.e. a curve of degree 8 in $\mathbb{P} (1,1,2)$). The
points of intersection of this curve with a generic ray is a
symmetric configuration of four points. Such configuration on a
generic ray does not have any point lying on the conic
$\widetilde{D}$ (as $R_8$ cuts $\widetilde{D}$ in 8 points) so do
not contain the center of symmetry (that is the vertex of the
cone). Thus the branch locus does not pass through the vertex of
the cone. After taking the quotient by the symmetry, we get
$\varphi$ to be a 2:1 morphism to $\mathbb{P}(1,1,4)$ with branch
locus which is a curve of degree 8 that does not pass through the
point (0:0:1). Let us consider the variety $\widetilde{Z}$ in
$\mathbb{P}(1,1,4,4)$ given by the equation as in case (a) from
the assertion, with $G$ defining the branch locus $R_8$. To prove
that we are in this case (a), we need to prove that $Z$ is
isomorphic to $\widetilde{Z}$. As $Z$ and $\widetilde{Z}$  have a
2:1 morphism onto $\mathbb{P}(1,1,4)$ with the same branch locus,
it is enough to prove that $\widetilde{Z}$ has only quotient
singularities, hence is normal. The only singularities of
$\widetilde{Z}$ are:

-The singularities of the branch locus - they are Du Val because
the branch locus has only simple singularities and do not pass
through the singular point $(0:0:1)$,

- and the points in the pre-image (by the projection) of the
singularity on the base - but the projection is a local (analytic)
isomorphism of their neighborhood with the singularity of
$\mathbb{P}(1,1,4)$ that is of type $\frac{1}{4} (1,1)$.

This shows that we are in case (a) of the theorem and completes
the proof of the first part.

\smallskip For the converse let us take a surface described in one
of the cases (a),  (b). We can check by adjunction formula that
they are del Pezzo surfaces (we have already checked that they are
normal). The only thing we need to prove is their type of
singularities (index $\leq 2$). We have computed that in case (a)
where the only non Du Val singularities are $\frac{1}{4} (1,1)$
that is of type $K_1$.

In case  (b) all singularities are Du Val except the singularity
in the point (0:0:0:1) that is the quotient by the symmetry $s$ of
the singularity of the quartic given by the equation
$\{F(a,b,c,d^2) =0\}$ in the point (0:0:0:1). To prove that it is
of type $K_n$, we need to prove that above quartic has a
$A_{2n+1}$ singularity in this point. This is a consequence of the
following lemma:
\begin{lem}\label{symmetricADE}
If $S$ is a surface in $\mathbb{P}^m$ with ADE singularities
symmetric with respect to a point $O$ (and a hyperplane) and
singular in this point, then it has a singularity $A_{2n+1}$ in
$O$.
\end{lem}
\proof[Proof of lemma] Observe that resolving the singularities of
$S$ (by blowing up step by step) we can extend the symmetry  to an
involution on the desingularization $\widetilde{S}$ of $S$. This
involution preserves the configuration of curves obtained as the
pre-image of our point $O$. Moreover the involution induced on
this configuration of curves has no isolated fixed points.
Moreover it preserves one of its curves (the first exceptional
divisor). We can see that from all ADE configurations only
$A_{2n+1}$ can have such an involution. Indeed:

\smallskip
 \textbf{If the configuration is of type $A_k$} then the involution
preserves each of the curves in the diagram. So each curve has at
least two fixed points on it. Moreover each curve is either fixed
(made of fixed points) or it intersects exactly two fixed curves,
and no two fixed curves intersect each other. Finally we get that
$k$ is odd.

\smallskip
\textbf{If the configuration is of type $D_k$ or $E_k$}, let us
consider the curve that cuts three other curves (let us call it
$L$). It has to be preserved by the involution. So it is fixed or
cuts at exactly two fixed curves. If it was fixed then, all curves
cutting it would be preserved by the involution and not fixed.
But, one of them does not intersect any other curves except $L$,
so it has to be fixed. This gives a contradiction. If it cuts two
fixed curves then the third curve cutting $L$ also has to be
preserved by the involution so it gives a third fixed point on
$L$. So $L$ is fixed and we get a contradiction.

This ends the proof of the lemma and in consequence the theorem.
\qed

\begin{rem}
All cases of above theorem belong to the same family of complete
intersections of hypersurfaces of degree 2 and 4 in the projective
space $\mathbb{P} (1,1,1,2,2)$.
\end{rem}

\begin{rem}\label{impossible}
We can observe that the hyperelliptic case with smooth image, from
the proof of the theorem, never occurs. This is because it would
be an index 1 case.
In this case $\psi \circ \varphi =\varphi_{|-2 K_Z|}$ (see the
proof of corollary \ref{weighted g3}). By using \cite{HW}, we then
conclude, that in the index 1 case with $g=3$ the system $|-2
K_Z|$ is very ample.

In other words, any del Pezzo surfaces of index 1 can also be
obtained as a double covering of $\mathbb{P}^2$ branched over a
quartic curve, but by our construction, it will be mapped directly
onto a quartic surface in $\mathbb{P}(1,1,1,2)$.
\end{rem}
\begin{rem}
Case (a) of Theorem \ref{g3} is the only example in the list of
del Pezzo surfaces from \cite{AN} (it has number 25 in Table 3),
with two index 2 singularities. Case  (b) of this theorem covers
all other cases of del Pezzo surfaces of index 2 with $g=3$.
Gorenstein surfaces are represented by quartics not passing
through the point $(0:0:0:1)$, while quartics passing through this
point have in this point their only index 2 singularity.
\end{rem}

On the other hand, Theorem \ref{g3} can be reformulated in the
following way:
\begin{cor}\label{weighted g3}
We have one of the following morphisms from $Z$ to a projective
space:
\begin{itemize}
 \item[(a)] a $2:1$ morphism onto the cone in $\mathbb{P}^5$, over
the rational normal quartic curve in $\mathbb{P} ^4$, branched
over the intersection of this cone with a quadric not passing
through its vertex,
 \item[ (b)] an isomorphism given by $|-2K_Z|$
to a complete intersection in $\mathbb{P}^6$ of a cone over the
Veronese surface with a (nonsingular) quadric.
\end{itemize}
\end{cor}

\begin{proof}
The corollary follows from Theorem \ref{g3} by composing our map
$\varphi$ from Lemma \ref{Key} with the embedding $\psi$.

From the proof of Lemma \ref{Key} we have that $( \psi \circ
\varphi \circ r) ^{\ast } ( \mathcal{O}_{\mathbb{P}^3}(1) )
\subset |C_3|$. Moreover, by Lemma \ref{dimensions of systems},
the system $|C_3|$ is of dimension 6 which is equal to the
dimension of the space spanned by the image $( \psi \circ \varphi
\circ r) (X)$ in case  (b). This means, that in this case $\psi
\circ \varphi= \varphi_{|-2K_Z|}$.

Case (a) from Theorem \ref{g3} leads directly to case (a) in above
corollary.
\end{proof}
\begin{rem}
From the above we deduce that in case (b) of Theorem \ref{g3}, the
linear system $|-2K_Z|$ is very ample. In case (a) this reasoning
doesn't work, nevertheless in this case we can also prove that
$|-2K_Z|$ is very ample. We prove that as follows. We can see that
$Z$ is a normal hypersurface of degree 8 in $\mathbb{P}(1,1,4,4)$.
Then by adjunction formula $|-2K_Z|=
\mathcal{O}_{\mathbb{P}(1,1,4,4)}(4)|_{Z}$. We know that the
latter is very ample.

\end{rem}
\end{section}

\begin{section}{Higher genus cases}\label{sec4}
In this section we work out the remaining cases, still using the
same method. The basic result is the following theorem about the
morphism given by the system $|-2K_Z |$ on $Z$.
\begin{thm}\label{isomorphic -2KZ}
If $Z$ is a log del Pezzo surface of index $ \leq 2 $ with
$K_{Z}^2 \geq 3$ then $|-2K_Z |$ gives an isomorphism onto the
image.
\end{thm}

Let us first prove some lemmas that will be important in the
further discussion.

\begin{lemm} \label{Scroll}
If $S$ is a smooth rational normal scroll of degree $d$ and
dimension $n$ in $\mathbb{P}^{n+d-1}$ invariant under a symmetry
with respect to a point $O$ (and some hyperplane), then $S$ is a
quadric hypersurface not passing through $O$.
\end{lemm}
\begin{proof}
Suppose that $S$ is a smooth scroll which is symmetric with
respect to a point $O$. By cutting $S$ with $n-d$ generic
hyperplanes passing through $O$, we reduce the problem to $S$
being a rational normal curve of degree $d$ in $\mathbb{P}^d$.

\smallskip If $S$ is a symmetric rational normal curve o degree $d$ in
$\mathbb{P}^d$, then a generic hyperplane $H$ passing through $O$
cuts $S$ in a symmetric configuration of $d$ points that span this
hyperplane. The only possibility for this to happen is that $d =2$
and $O$ does not lie on $S$.
\end{proof}

\begin{lemm}\label{sym}
The ideal of a symmetric variety $V$ is generated by symmetric
elements.
\end{lemm}
\begin{proof}
For each $n$ let us consider the linear system of hypersurfaces of
degree $n$ containing $V$. This is a projective space with an
involution on it (given by the symmetry preserving $V$). But an
involution on a projective space is a symmetry. To finish the
proof it is enough to observe that the fixed points of any
symmetry on a projective space generate this space.
\end{proof}

Let us pass to the proof of the theorem.
\begin{proof}[Proof of Theorem \ref{isomorphic -2KZ}]
Using Lemma \ref{Key}, to prove the first statement it is enough
to check that the map $\varphi_{|D_g|}$ is birational and deduce
that $\varphi$ is an isomorphism onto its image. Next by comparing
dimensions as in preceding theorems we obtain, that $\varphi$ is
given by the system $|-2K_Z|$.

Suppose that $|D_g|$ is hyperelliptic. Since $D_g$ have genus $>3$
according to \cite[th.5.7]{SD} we have three possibilities for the
morphism given by this system:

\smallskip
$\bullet$ It is a double covering of a smooth centrally symmetric
rational normal scroll of degree $g-1$ in $\mathbb{P}^g$. This is
impossible because smooth scrolls of codimension 2 are never
centrally symmetric (see Lemma \ref{Scroll})

\smallskip
$\bullet$ It is a double covering of a cone over a twisted cubic.
In this case $g=4$ and $|D_4|=|3E +2\Gamma _1 +\Gamma _2 |$ where
$E$ is an elliptic curve and $\Gamma _1$, $\Gamma _2$ are smooth
and rational curve such that $E. \Gamma _1 =\Gamma _1 . \Gamma _2
=1$. The involution on the K3 surface $X$ has to preserve the
curves $\Gamma _1$, $\Gamma _2$ (because it preserves the system,
the elliptic curve $E$ and the inverse image of the vertex of the
cone) and restricted to them cannot have isolated fixed points (as
$D_g . \Gamma _i =0 $ for $i=1,2$). This is impossible because on
each of these two curves there have to be at least two fixed
points, but two fixed curves cannot intersect.

\smallskip
$\bullet$  It is a double covering of a cone over a twisted
quartic. In this case $g=5$ and $|D_5|=|4E +2\Gamma|$ where $E$ is
an elliptic curve and $\Gamma $ is a smooth rational curve such
that $E. \Gamma  =1$. Then the map given by $D_5$ can be extended
to a map $\widetilde{\varphi}$ on the desingularization of the
cone, that is the scroll $\mathbb{F}_4$. The map
$\widetilde{\varphi}$ has to be branched over the exceptional
section $e$ of the scroll (as its inverse image on the K3 surface
$X$ is an irreducible $-2$ curve). Moreover we know that the
branch locus of this map has to cut the generic line from the
ruling of the scroll in a configuration of 4 points symmetric with
respect to the 2 points of intersection of this line with the
section $e$ and the image of $D_5$. It also cuts $D_5$ in 12
points. These statements together give us a contradiction because
the generic line of the ruling does not pass through any of these
12 points, so the configuration of the three remaining points on
the generic line cannot be symmetric.
\end{proof}

We can know pass to a deeper analysis of particular cases
considered in this section.
\begin{thm}\label{g4}\ \\
\textbf{\emph{1)}} If $Z$ is a log del Pezzo of index $ \leq 2 $
with $g=K_{Z}^2 +1= 4$, it is isomorphic to a complete
intersection in w.p.s. $\mathbb{P}(1,1,1,1,2) \subset
\mathbb{P}^{10} $ (with coordinates $x,y,z,t,u$) defined by
polynomials $F, G$ of weighted degree 2 and 3 respectively, such
that the variety in $\mathbb{P}^4$ (with coordinates $a,b,c,d,e$)
given by equations $F(a,b,c,d,e^2 )=0$, $G(a,b,c,d,e^2 )=0$  has
only Du Val
singularities lying outside the hyperplane $\{e=0\}$.\\
 \textbf{\emph{2)}} Conversely every such complete
intersection is a log del Pezzo surface of index $\leq 2$.
\end{thm}

\begin{proof}
After the above discussion the proof of Theorem \ref{g4} is
straightforward. It is enough to observe (using \cite[6.5.3]{SD})
that the morphism given by $D_4$ is birational onto the complete
intersection of a symmetric quadric and a symmetric cubic in
$\mathbb{P}^4$. We then take the quotient by this symmetry and
finish the proof.
\end{proof}

\begin{rem}
From Lemma \ref{Scroll} we deduce that the log del Pezzo surface
obtained in Theorem \ref{g4} has an index 2 singularity if and
only if the quadric given by $F(a,b,c,d,e^2 )=0$ is singular. If
this quadric is smooth, then $Z$ is isomorphic to a cubic surface
in $\mathbb{P}^3$.
\end{rem}
Recall (see \cite{HW}) that, if $Z$ is of index 1 (Gorenstein) and
$g\geq5$, then it is isomorphic to an intersection of
$\frac{(g-1)(g-4)}{2}$ hypersurfaces of degree 2 in
$\mathbb{P}^{g-1}$. In the next theorem we deal with the cases of
index $=2$.
\begin{thm}\label{higher genus}
If $Z$ is a log del Pezzo surface of index $=2$ and with $K_{Z}^2
\geq 4$, then the image $\varphi_{|-2K_Z|} (Z)$:
 \begin{itemize}
 \item[(a)]is a cone over a rational normal quartic in
 $\mathbb{P}^5$,
 \item[(b)] is contained in a threefold cone over the image by the double
Veronese embedding of a rational normal scrollar surface,
 \item[(c)] is contained in a threefold cone over the image of the quadruple
Veronese embedding of $\mathbb{P}^2$ in $\mathbb{P}^{14}$.
\end{itemize}
\end{thm}
\begin{proof}
Observe that by \cite[th.7.2]{SD} we know that the ideal $I$ of
the image $\varphi_{|D_g|} (X)$ is generated by elements of degree
2 and 3. Since the index 2 singularity can appear only if
$\varphi_{|D_g|} (X)$ passes through the center of symmetry $O$
all quadrics and cubics generating the ideal have to pass through
$O$ and by Lemma \ref{symmetry of the image} have to define
symmetric varieties. The only possibility that a quadric $Q$ is
symmetric with respect to a point lying on it is that $Q$ is a
cone. Following \cite{SD} let us denote the variety given by all
elements of degree 2 in $I$ by $T_{|D_g|}$. Naturally
$\varphi_{|D_g|} (X) \subset T_{|D_g|}$. According to \cite{SD}
and taking into account that $T_{|D_g|}$ has to be symmetric we
have two possibilities:
\begin{itemize}
 \item either $\varphi_{|D_g|} (X) = T_{|D_g|}$
 \item or $T_{|D_g|}$ is a threefold cone
\end{itemize}

\textbf{In the first case} $\varphi_{|D_g|} (X)$ has to be a cone
with Du Val singularities (in particular isolated). Hence, it has
only one singularity, that is the one in the vertex of the cone.
It has to be a singularity $A_1$ since the exceptional curve of
the blowing up of an irreducible cone is irreducible. Therefore
$Z$ is a cone with a singularity $K_1$ in the vertex. Using
\cite[table 3]{AN} we get that the DPN surface $Y$ is isomorphic
to $\mathbb{F}_4$. By contracting the exceptional section we
obtain $Z$ to be the cone over a rational normal quartic as in
case (a) of the theorem.

\textbf{In the second case} following \cite{SD} and taking into
 account that $T_{|D_g|}$ is a cone we obtain that
 $T_{|D_g|}$ has to be one of the following:
 \begin{itemize}
 \item[-] a cone over a rational normal scrollar surface
 \item[-] a cone over the Veronese surface
 \end{itemize}
This proves the theorem since after taking the quotient by the
symmetry above cases pass to cases (b) and (c) from the assertion
respectively.
\end{proof}
As a consequence of the proof of above theorem we get the
following corollary.
\begin{cor} If $Z$ is a log del Pezzo surface of index $=2$ with
$g=K_Z^2+1\geq4$, then it is isomorphic to a surface in
$\mathbb{P}(1^g,2)$ generated by equations of degree 2 and 3.
\end{cor}

In the following theorem we write case (c) from Theorem
\ref{higher genus} in terms of equations in a weighted projective
space in a different and more precise way.

\begin{thm}\ \\
\textbf{\emph{1)}} If $Z$ is a log del Pezzo surface corresponding
to case (c) of the Theorem \ref{higher genus}, then $Z$ is
isomorphic to a surface in $\mathbb{P}(1,1,1,4)$ given by an
equation $F_5$ such that $F_5 (a,b,c,d^4)=0$ defines a surface of
degree 5 in $\mathbb{P}^3$ (with coordinates $a$, $b$, $c$, $d$)
with ADE singularities and smooth in the point $(0:0:0:1)$.\\
 \textbf{\emph{2)}} The converse statement is also true.
\end{thm}
\begin{proof} In this case, the morphism $\varphi_{|D_6|} \colon X \longrightarrow
\mathbb{P}^6$ is birational onto its image, that is contained in a
cone over the Veronese surface. This means that
$\varphi_{|D_6|}(X)$ is a K3 surface with Du Val singularities
lying on a cone over the Veronese surface (that is in
$\mathbb{P}(1,1,1,2)$). Hence it is a hypersurface of degree 5 in
$\mathbb{P}(1,1,1,2)$. Let us denote its equation by
$\widetilde{F}_5$.
Observe that the surface $S$, given in $\mathbb{P}^3$ (with
coordinates $a$, $b$, $c$, $d$) by the equation
$\widetilde{F}_5(a,b,c,d^2)=0$, is normal as it is a hypersurface
with isolated singularities. The singularity of this surface in
the point $(0:0:0:1)$ is thus a normal double covering of a Du Val
singularity branched in one point, hence it is a Du Val
singularity (or it can be smooth). From Lemma \ref{symmetricADE}
we deduce that $S$ is either smooth in $(0:0:0:1)$ or has there an
$A_{2n+1}$ singularity. The latter is impossible, since the
quotient of an $A_{2n+1}$ singularity by a central symmetry is of
type $K_n$ (that is not ADE). We deduce that $S$ is smooth in the
point $(0:0:0:1)$ and has only Du Val singularities outside.

We know from \ref{symmetry of the image} that $\varphi_{|D_6|}(X)$
is irreducible and symmetric with respect to the symmetry
$\widetilde{s}\colon \mathbb{P}(1,1,1,2)\ni (x:y:z:t)\mapsto
(x:y:z:-t) \in \mathbb{P}(1,1,1,2)$. In the same way as in
theorems \ref{M} and \ref{g3} this induces that $\widetilde{F}_5$
is of the form $F_5(x,y,z,t^2)$. By taking the quotient by the
symmetry $s$, we obtain the first part of the theorem.

\smallskip
The converse statement is a consequence of the adjunction formula
and the normality of above equations.
\end{proof}
\begin{rem} From the above theorem we deduce, that all surfaces
corresponding to case $(c)$ of Theorem \ref{higher genus} have a
$K_1$ singularity.
\end{rem}
\begin{rem}
To distinguish between cases (a), (b) and (c) from Theorem
\ref{higher genus}, we can use \cite{SD} once more.

Case (a) is well described and isolated. It is the only case where
$g=10$. In table 3 of \cite{AN} it has number 11

Case (b) corresponds to the case where there exists an elliptic
pencil $|E|$ on the K3 surface $X$ such that $E.D_g =3$ (it means
that we are in the trigonal case, i. e. the curve $D_g$ is
trigonal).

Case (c) corresponds to the case where $D_g$ is isomorphic to a
plane quintic. In this case $g=6$ and $|D_6|=|2B + \Gamma |$ where
$B$ is some curve of genus 2 and $\Gamma$ a rational curve such
that $B.\Gamma =1$.

Let us see which of the cases from \cite[table 3]{AN} corresponds
to case (b) and  (c). To do this we check all diagrams with k>0
and $g=6$ together with all diagrams that we obtain by choosing
subgraphs of above ones. We see that, except in case 15, in all
these cases there is always one of the configuration of curves
written below.

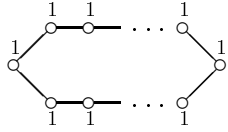
\begin{itemize}
 \item[(i)]$$
 \setlength{\unitlength}{0.5mm}
 \begin{picture}(70,30)(0,-20)
 \put(0,0){\circle{3}}
 \put(1,1){\line(1,1){8}}
 \put(10,10){\circle{3}}
 \put(20,10){\circle{3}}
 \put(45,10){\circle{3}}
 \put(46,9){\line(1,-1){8}}
 \put(55,0){\circle{3}}
 \put(11.5,10){\line(1,0){7}}
 \put(21.5,10){\line(1,0){7}}
\put(1,-1){\line(1,-1){8}}
 \put(10,-10){\circle{3}}
 \put(20,-10){\circle{3}}
 \put(45,-10){\circle{3}}

 \put(-1,3){\makebox{\tiny{1}}}
 \put(9,13){\makebox{\tiny{1}}}
 \put(19,13){\makebox{\tiny{1}}}
 \put(44,13){\makebox{\tiny{1}}}
 \put(54,3){\makebox{\tiny{1}}}
 \put(9,-16){\makebox{\tiny{1}}}
 \put(19,-16){\makebox{\tiny{1}}}
 \put(44,-16){\makebox{\tiny{1}}}

 \put(46,-9){\line(1,1){8}}
 \put(11.5,-10){\line(1,0){7}}
 \put(21.5,-10){\line(1,0){7}}
 \put(36.5,9){\makebox(0,0)[b]{\dots}}
 \put(36.5,-11){\makebox(0,0)[b]{\dots}}
\end{picture}$$
 \item[(ii)]$$
\setlength{\unitlength}{0.5mm}
 
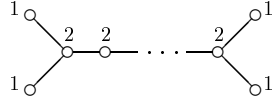
\begin{picture}(50,10)(0,0)
 \put(0,20){\circle{3}}
 \put(0,0){\circle{3}}
 \put(10.0,10.0){\circle{3}}
 \put(20.0,10.0){\circle{3}}
 \put(50.0,10.0){\circle{3}}
 \put(60.0,20.0){\circle{3}}
 \put(60.0,0){\circle{3}}
 \put(1,1){\line(1,1){8}}
 \put(1,19){\line(1,-1){8}}
 \put(11.5,10){\line(1,0){7}}
 \put(21.5,10){\line(1,0){7.0}}
 \put(41.5,10){\line(1,0){7.0}}
 \put(51,11){\line(1,1){8}}
 \put(51,9){\line(1,-1){8}}

 \put(-5.5,20){\makebox{\tiny{1}}}
 \put(-5.5,0){\makebox{\tiny{1}}}
 \put(9,13.0){\makebox{\tiny{2}}}
 \put(19.0,13.0){\makebox{\tiny{2}}}
 \put(49.0,13.0){\makebox{\tiny{2}}}
 \put(62,20.0){\makebox{\tiny{1}}}
 \put(62,0){\makebox{\tiny{1}}}

 \put(36,10){\makebox(0,0){\dots}}
\end{picture}$$
\item[(iii)]
$$\setlength{\unitlength}{0.5mm}
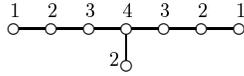
\begin{picture}(70,25)(0,0)
 \put(0,20){\circle{3}}
 \put(10,20){\circle{3}}
 \put(20,20){\circle{3}}
 \put(30,20){\circle{3}}
 \put(40,20){\circle{3}}
 \put(50,20){\circle{3}}
 \put(60,20){\circle{3}}
 \put(30,10){\circle{3}}

 \put(-1,23){\makebox{\tiny{1}}}
 \put(9,23){\makebox{\tiny{2}}}
 \put(19,23){\makebox{\tiny{3}}}
 \put(29,23){\makebox{\tiny{4}}}
 \put(39,23){\makebox{\tiny{3}}}
 \put(49,23){\makebox{\tiny{2}}}
 \put(59,23){\makebox{\tiny{1}}}
 \put(25.5,10){\makebox{\tiny{2}}}

 \put(1.5,20){\line(1,0){7}}
 \put(11.5,20){\line(1,0){7}}
 \put(21.5,20){\line(1,0){7}}
 \put(31.5,20){\line(1,0){7}}
 \put(41.5,20){\line(1,0){7}}
 \put(51.5,20){\line(1,0){7}}
 \put(30,18.5){\line(0,-1){7}}
 \end{picture}
$$
\item[(iv)]\setlength{\unitlength}{0.5mm}
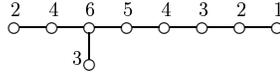
$$\begin{picture}(70,25)(0,0)
 \put(10,20){\circle{3}}
 \put(20,20){\circle{3}}
 \put(30,20){\circle{3}}
 \put(40,20){\circle{3}}
 \put(50,20){\circle{3}}
 \put(60,20){\circle{3}}
 \put(70,20){\circle{3}}
 \put(80,20){\circle{3}}

 \put(30,10){\circle{3}}

 \put(9,23){\makebox{\tiny{2}}}
 \put(19,23){\makebox{\tiny{4}}}
 \put(29,23){\makebox{\tiny{6}}}
 \put(39,23){\makebox{\tiny{5}}}
 \put(49,23){\makebox{\tiny{4}}}
 \put(59,23){\makebox{\tiny{3}}}
 \put(69,23){\makebox{\tiny{2}}}
 \put(79,23){\makebox{\tiny{1}}}

 \put(25.5,10){\makebox{\tiny{3}}}

 \put(11.5,20){\line(1,0){7}}
 \put(21.5,20){\line(1,0){7}}
 \put(31.5,20){\line(1,0){7}}
 \put(41.5,20){\line(1,0){7}}
 \put(51.5,20){\line(1,0){7}}
 \put(61.5,20){\line(1,0){7}}
 \put(71.5,20){\line(1,0){7}}

 \put(30,18.5){\line(0,-1){7}}
 \end{picture}
$$
\item[(v)]\setlength{\unitlength}{0.5mm}
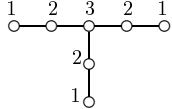
$$\begin{picture}(70,25)(0,0)
 \put(10,20){\circle{3}}
 \put(20,20){\circle{3}}
 \put(30,20){\circle{3}}
 \put(40,20){\circle{3}}
 \put(50,20){\circle{3}}
 \put(30,10){\circle{3}}
 \put(30,0){\circle{3}}
 \put(8,23){\makebox{\tiny{1}}}
 \put(19,23){\makebox{\tiny{2}}}
 \put(29,23){\makebox{\tiny{3}}}
 \put(39,23){\makebox{\tiny{2}}}
 \put(48,23){\makebox{\tiny{1}}}
 \put(25.5,10){\makebox{\tiny{2}}}
 \put(25,0){\makebox{\tiny{1}}}
 \put(11.5,20){\line(1,0){7}}
 \put(21.5,20){\line(1,0){7}}
 \put(31.5,20){\line(1,0){7}}
 \put(41.5,20){\line(1,0){7}}
 \put(30,18.5){\line(0,-1){7}}
 \put(30,8.5){\line(0,-1){7}}
 \end{picture}
$$
\end{itemize}

We wrote those configurations of curves with some multiplicities,
so that they generate an appropriate linear system. We check that
the linear systems of each of these configurations of curves with
chosen multiplicities is an elliptic pencil on X (we do it by
intersecting it with each of its components and getting always
zero). In all above cases, we can choose these configurations so
that $E.C_g =3$.

In case 15 for each of the diagrams (there is one main diagram but
there are also diagrams generated by subgraphs), using the method
described at the end of section \ref{intro}, we can write the
system $|D_g|$ as $|2B + \Gamma|$. For example in the case
represented by the full diagram we write $|D_g|$ in the following
way:

 \vskip10pt
 \setlength{\unitlength}{0.5mm}
 $\begin{picture}(230,25)(20,0)
 \put(20,10){\circle{3}}
 \put(30,10){\circle{3}}
 \put(30,20){\circle{3}}
 \put(30,0){\circle{3}}
 \put(40,20){\circle{3}}
 \put(50,20){\circle{3}}
 \put(60,20){\circle{3}}
 \put(40,0){\circle{3}}
 \put(50,0){\circle{3}}
 \put(60,0){\circle{3}}
 \put(21.5,10){\line(1,0){7}}
 \put(31.5,20){\line(1,0){7}}
 \put(41.5,20){\line(1,0){7}}
 \put(51.5,20){\line(1,0){7}}
 \put(31.5,0){\line(1,0){7}}
 \put(41.5,0){\line(1,0){7}}
 \put(51.5,0){\line(1,0){7}}
 \put(30,1.5){\line(0,1){7}}
 \put(30,18.5){\line(0,-1){7}}
 \put(19,3){\makebox{\small{5}}}
 \put(33,9){\makebox{\small{10}}}
 \put(30,23){\makebox{\small{8}}}
 \put(30,-7){\makebox{\small{8}}}
 \put(40,23){\makebox{\small{6}}}
 \put(50,23){\makebox{\small{4}}}
 \put(60,23){\makebox{\small{2}}}
 \put(40,-7){\makebox{\small{6}}}
 \put(50,-7){\makebox{\small{4}}}
 \put(60,-7){\makebox{\small{2}}}
 \put(90,10){\makebox(0,0)[c]{\huge{=}}}
 \put(130,10){\makebox(0,0)[c]{\huge{+}}}
 \put(150,10){\makebox(0,0)[c]{\huge{2}}}
 \put(170,10){\makebox(0,0)[c]{\huge{$\times$}}}
 \put(110,10){\circle{3}}
 \put(190,10){\circle{3}}
 \put(200,10){\circle{3}}
 \put(200,20){\circle{3}}
 \put(200,0){\circle{3}}
 \put(210,0){\circle{3}}
 \put(210,20){\circle{3}}
 \put(220,20){\circle{3}}
 \put(230,20){\circle{3}}
 \put(230,0){\circle{3}}
 \put(220,0){\circle{3}}
 \put(230,0){\circle{3}}
 \put(191.5,10){\line(1,0){7}}
 \put(201.5,20){\line(1,0){7}}
 \put(211.5,20){\line(1,0){7}}
 \put(221.5,20){\line(1,0){7}}
 \put(201.5,0){\line(1,0){7}}
 \put(211.5,0){\line(1,0){7}}
 \put(221.5,0){\line(1,0){7}}
 \put(200,1.5){\line(0,1){7}}
 \put(200,18.5){\line(0,-1){7}}
 \put(190,3){\makebox{\small{2}}}
 \put(203,9){\makebox{\small{5}}}
 \put(200,23){\makebox{\small{4}}}
 \put(200,-7){\makebox{\small{4}}}
 \put(210,23){\makebox{\small{3}}}
 \put(220,23){\makebox{\small{2}}}
 \put(230,23){\makebox{\small{1}}}
 \put(210,-7){\makebox{\small{3}}}
 \put(220,-7){\makebox{\small{2}}}
 \put(230,-7){\makebox{\small{1}}}
 \end{picture}$

\vskip20pt

\end{rem}


\begin{rem}
Let $Z$ be a log del Pezzo surface of index $=2$ and with
$g=K_{Z}^2 +1=5$. Then the image of the map $|D_5 |$ is a degree
eight intersection in $\mathbb{P}^5$ of a cone over a rational
normal scrollar surface in $\mathbb{P}^4$ (it is of degree 3) with
two symmetric (see Lemma \ref{sym}) cubic hypersurfaces (in fact
it is a component of the intersection of this cone with one cubic
hypersurface where the second component has to be a plane). We can
also obtain a converse statement saying that: each surface with Du
Val singularities and obtained as a component of degree 8 of a
section of the cone with a cubic surface is a K3 surface
(singular). We add to this the condition to being symmetric. Then
if we take the quotient by the symmetry we obtain a description of
del Pezzo surfaces of index 2 and $K_{Z}^2 =4$.
\end{rem}

In general if $Z$ has $g \geq 6$ our method does not give
satisfactory results because in those cases we do not have a list
of corresponding $K3$ surfaces given in terms of equations. For
this reason we give the complete list of del Pezzo surfaces with
$g\geq 6$ (see \cite[table 3]{AN}). We do the following.

\begin{thm}
There are only five families of dimension $\geq 1$ of del Pezzo
surfaces with index 2 and $g \geq 6$ and a finite number of
isomorphism classes not in these families.
\end{thm}
\proof We use the toric argument described in \cite[proof of thm
4.2]{AN} to all the diagrams of DPN (see \cite[table 3]{AN})
surfaces with $g\geq 6$. We do as follows. We consider all
diagrams of exceptional curves that occur on DPN surfaces with
$g\geq 6$, that is all diagrams from \cite[table 3]{AN} with
$g\geq 3$ and all their associated diagrams. The diagram permits
us to write the DPN surfaces as obtained by a sequence of blowing
ups of $\mathbb{P}^2$, $\mathbb{P}^1 \times \mathbb{P}^1$ or
$\mathbb{F}_n$. These are toric varieties. The torus group acts on
the space of choices for the blowing ups. Our task is to compute
the dimensions of spaces of orbits.  After a straightforward
computation, we obtain that the only non discrete families of DPN
surfaces of $g\geq 6$ with a fixed diagram of exceptional curves
are:
\begin{itemize}
\item[(1)] $\mathbb{P}^2$ blown up in five point on a line.

\item[(2)] $\mathbb{F}_1$ blown up in four points on a fiber of
which one on the exceptional section.

\item[(3)] $\mathbb{F}_3$ blown up in four points on a fiber of
which one on the exceptional section.

\item[(4)]$\mathbb{P}^1 \times \mathbb{P}^1$ blown up in four
points on a fiber.

\item[(5)] $\mathbb{F}_2$ blown up in four points on a fiber of
which one on the exceptional section and another point on this
section.
\end{itemize}
The first of them is two dimensional, because it is classified by
configurations of five point on a line. The others are one
dimensional, as they are connected with configurations of four
points on a line. By contracting $\mathrm{Duv}$ and $\mathrm{Log}$
parts of the diagrams we obtain families of log del Pezzo
surfaces. \qed
\end{section}
\begin{rem} The families (1) and (2) from above theorem come from
subgraphs from the case 15 of \cite[table 3]{AN}, so their
corresponding families of del Pezzo surfaces are described as
quintics in $\mathbb{P}(1,1,1,4)$.
\end{rem}

\begin{section}{Application}\label{sec5}
As an application of these results with the help of \cite[table
3]{AN} we can get the list of all configurations of singularities
that can occur on surfaces defined by equations described in the
paper, or equivalently on symmetric singular K3 surfaces
corresponding to them. Moreover from Section \ref{intro} we know
that all these configurations of singularities do occur. Below we
show some interesting examples of this type.
\begin{ex} As an example we can get a list
of all configurations of singularities that appear on a symmetric
sextic plane curve with at most a double singularity in the center
of symmetry. This result without the assumption of being symmetric
is described in \cite{AN}. The complete list described there is
much longer (and not formulated explicitly). We first deal with
the case where the sextic is symmetric with respect to a point on
it. We have the following list:
\begin{itemize}
\item[-]An $A_1$ singularity in the center of symmetry and twice
(because singular points appear in pairs) the following
configuration (compare with \cite[sec.4.2]{AN}):\\
$2D_4$;\\
$D_8$;\\
$D_7$;\\
$D_6 2A_1$, $D_6 A_1$, $D_6$; $D_5 A_3$, $D_5 A_2$, $D_5 2A_1$,
$D_5
A_1$, $D_5$;\\
$D_4 A_3$, $D_4 A_2$, $D_4 3A_1$, $D_4 2A_1$, $D_4 A_1$, $D_4$;\\
$A_7$;\\
$A_6$;\\
$A 5 2A_1$, $A_5 A_1$, $A_5$;\\
$A_4 A_3$, $A_4 A_2$, $A_4 2A_1$, $A_4 A_1$, $A_4$;\\
$2A_3 A_1$, $2A_3$, $A_3 A_2 2A_1$, $A_3 A_2 A_1$, $A_3 A_2$, $A_3
4A_1$, $A_3 3A_1$, $A_3 2A_1$, $A_3 A_1$, $A_3$;\\
$2A_2 2A_1$, $2A_2 A_1$, $2A_2$, $A_2 4A_1$, $A_2 3A_1$, $A_2
2A_1$, $A_2 A_1$, $A_2$;\\
$6A_1$, $5A_1$, $4A_1$, $3A_1$;\\
$2A_1$, $A_1$;\\
$\emptyset$;
 \item[-] A singularity $A_3$ in the center of symmetry
and twice the following configuration:\\
$A_7$;\\
$A_6$;\\
$A_5 A_1$, $A_5$;\\
$A_4 A_2$, $A_4 A_1$, $A_4$;\\
$2A_3$, $A_3 A_2$, $A_3 2A_1$, $A_3 A_1$, $A_3$;\\
$2A_2 A_1$, $2A_2$, $A_2 2A_1$, $A_2 A_1$, $A_2$;\\
$4A_1$, $3A_1$, $2A_1$, $A_1$;\\
$\emptyset$.
 \item[-] A singularity $A_5$ in the center of symmetry
and twice the following configuration:\\
$A_5 A_1$, $A_5$;\\
$A_4 A_1$, $A_4$;\\
$A_3 2A_1$, $A_3 A_1$, $A_3$;\\
$2A_2 A_1$, $2A_2$, $A_2 2A_1$, $A_2 A_1$, $A_2$;\\
$4A_1$, $3A_1$, $2A_1$, $A_1$;\\
$\emptyset$.
 \item[-] A singularity $A_7$ or $A_9$ in the center of symmetry
and twice the following configuration:\\
$A_4$;\\
$A_3$;\\
$A_2 A_1$, $A_2$;\\
$2A_1$, $A_1$;\\
$\emptyset$;
 \item[-] A singularity $A_{11}$ in the center of symmetry
and twice the following configuration:\\
$A_2$;\\
$A_1$;\\
$\emptyset$;
 \item[-] A singularity $A_{13}$ in the center of symmetry
and twice $A_1$ or $\emptyset$;
 \item[-] A singularity $A_{15}$ or $A_{17}$ in the center of symmetry
and smooth outside.

\end{itemize}
If the sextic does not pass through the center of symmetry it is
either the sum of three tangent conics or has only simple
singularities. In the latter we have a very long list of cases,
described by all possible subgraphs of the
following graphs:\\
$E_8$, $A_8$, $A_7 A_1$, $A_5 A_2 A_1$, $2A_4$, $D_8$, $D_5 A_3$,
$E_6 E_2$, $E_7 A_1$, $D_6 2A_1$, $2D_4$, $2A_3 2A_1$, $4A_2$.

\end{ex}
\end{section}
\begin{ex}
As in the previous example we get a description of all
configurations of only Du Val singularities that appear on a
symmetric quartic in $\mathbb{P}^3$.
\begin{itemize}
 \item[-]
if the quartic passes through the center of symmetry there can be
only $A_n$ singularities and only in the following
configurations:\hfill \\
 $A_1$ in the center and a symmetric configuration of: \hfill \\
 $\emptyset $, $2A_7$, $2A_6$, $2A_5$, $2A_4$, $2A_3$, $2A_2$, $2A_1$,\\
 $2A_5 2A_1$, $2A_42A_2$, $4A_3$, $2A_4 2A_1$, $2A_3 2A_1$, $2A_2
 2A_1$, $4A_1$, $2A_3 2A_2$, $4A_2$,\\
 $2A_3 4A_1$, $2A_2 4A_1$, $6A_1$, $4A_2 2A_1$,\\
 $8A_1$,\\
 or $A_3$ in the center and a symmetric configuration of:\\
 $\emptyset $, $2A_3$, $2A_2$, $2A_1$,\\
 $4A_3$, $2A_3 2A_2$, $2A_3 2A_1$ $4A_2$, $2A_2 2A_1$, $4A_1$,\\
 $2A_3 4A_1$, $2A_2 4A_1$, $6A_1$,\\
 $8A_1$,\\
 or $A_5$ in the center and a symmetric configuration of:\\
 $\emptyset $, $2A_2$, $2A_1$,\\
 $4A_2$, $2A_2 2A_1$, $4A_1$,\\
 or $A_7$ in the center and a symmetric configuration of:\\
 $\emptyset $, $2A_2$, $2A_1$,\\
  $2A_2 2A_1$, $4A_1$,\\
 or $A_9$ in the center and a symmetric configuration of:\\
 $\emptyset $, $2A_2$, $2A_1$,\\
 or $A_11$ in the center and a symmetric configuration of:\\
 $\emptyset $, $2A_2$, $2A_1$\\
 \item[-] if the quartic does not pass through the center of
symmetry the configurations are described in \cite{HW} and
correspond to subgraphs of the extended Dynkin diagram
$\widetilde{E_7}$. We can also prove this looking to \cite[table
3]{AN}.
\end{itemize}
\end{ex}
\begin{ex} We know that the case (a) of Theorem \ref{g3} is the only
case with two index 2 singularities, and this is case number 25
from \cite[table 3]{AN}. The Du Val singularities of these log del
Pezzo surfaces correspond to the singularities of the branch
divisor of the map $\varphi$. From this we deduce, that all
configurations of singularities, that appear on a curve of degree
8 in $\mathbb{P}(1,1,4)$ having only Du Val singularities and not
passing through the point $(0:0:1)$, are defined by subgraphs
of the singularity $A_7$. That is, are one of the following: \\
$A_7$, $A_6$,\\
$A_5 A_1$, $A_4 A_1$, $A_3 A_1$, $A_2 A_1$, $2 A_1$, $A_4 A_2$,
$A_3 A_2$, $A_2 A_2$, $A_3 A_3$,\\
$A_3 2A_1$, $A_2 2A_1$, $2A_2 A_1$, $3A_1$,\\
$4A_1$.
\end{ex}
\begin{ex}
We can write down all configurations of only Du Val singularities
that occur on a symmetric surface being a complete intersection of
a quadric and a cubic in $\mathbb{P}^4$. They can have only $A_n$
singularities and only in the following configurations:
\begin{itemize}
 \item[-] if the surface passes through the center of symmetry
 $A_1$ in the center and a symmetric configuration of:\hfill \\
 $\emptyset $, $2A_5$, $2A_4$, $2A_3$, $2A_2$, $2A_1$,\hfill \\
 $2A_5 2A_1$, $2A_4 2A_1$, $2A_3 2A_1$, $2A_2
 2A_1$, $4A_1$, $4A_2$,\\
 $2A_3 4A_1$, $2A_2 4A_1$, $6A_1$, $4A_2 2A_1$,\\
 $8A_1$,\\
 or $A_3$ in the center and a symmetric configuration of:\\
 $\emptyset $, $2A_2$, $2A_1$,\\
 $4A_2$, $2A_2 2A_1$, $4A_1$\\
 or $A_5$ in the center and a symmetric configuration of:\\
 $\emptyset $, $2A_2$, $2A_1$,\\
 $4A_2$, $2A_2 2A_1$, $4A_1$,\\
 or $A_7$ in the center and a symmetric configuration of:\\
 $\emptyset $, $2A_1$,\\
 or $A_9$ in the center and a symmetric configuration of:\\
 $\emptyset $, $2A_1$,\\

 \item[-]
If the surface does not pass through the center of symmetry,  then
its quotient by the symmetry $s$ is a double covering of a cubic
surface in $\mathbb{P}^3$ branched over a smooth curve. This
implies that configurations of singularities that can appear on
this surface are non-branched coverings of the configurations of
singularities that can appear on the cubic. These are symmetric
configurations of subgraphs of the following graphs:\\
$E_6$, $A_5 A_1$, $3A_2$.

\end{itemize}
\end{ex}
\begin{ex}
The configurations of singularities that occur on a quintic with
only Du Val singularities outside the point $(0:0:0:1)$ in
$\mathbb{P}(1,1,1,4)$ and quasi-smooth in this point are:
\begin{itemize}\item[] A singularity $K_1$ in the point $(0:0:0:1)$ and
one of the following configurations of $A_n$ singularities:\\
$\emptyset$, $A_4$, $A_3$, $A_2$, $A_1$,\\
$A_2 A_1$, $2A_1$

\end{itemize}
\end{ex}

\vskip10pt
\begin{minipage}{8cm} Grzegorz Kapustka \\
Jagiellonian University\\
ul. Reymonta 4\\
30-059 Kraków\\
Poland\\
email: Grzegorz.Kapustka@im.uj.edu.pl
\end{minipage}
\begin{minipage}{7cm}
 Michał Kapustka\\
 Jagiellonian University\\
 ul. Reymonta 4\\
 30-059  Kraków\\
 Poland\\
 email: Michal.Kapustka@im.uj.edu.pl
\end{minipage}


\begin{thebibliography}{AN2}


 \bibitem[AN]{AN}
V. Alekseev, V. Nikulin, Classification of log del Pezzo surfaces
of index $\leq 2$, Preprint (2004), math.AG/0406536, 141 pages.
 \bibitem[AN2]{AN2}
V. Alekseev, V. Nikulin, Classification of del Pezzo surfaces with
log-terminal singularities of index = 2 and involutions of K3
surfaces, Soviet Math. Dokl. 39 (1989), no. 3, 507–511.
 \bibitem[D]{D}
M. Demazure, Surfaces de del Pezzo. II — IV, In M. Demazure, H.
Ch. Pinkham, and B. Teissier (eds.), S\'eminaire sur les
Singularit\'es des Surfaces, Lecture Notes in Mathematics, vol.
777, Springer, Berlin, 1980, Held at the Centre de Math\'ematiques
de l’  \'ecole Polytechnique, Palaiseau, 1976–1977, 23-69.
 \bibitem[DVa]{DVa}
P. Du Val, On isolated singularities of surfaces which do not
affect the conditions of adjunction. II, III, Proc. Cambridge
Phil. Soc. 30 (1934), 453-465.
 \bibitem[DVb]{DVb}
P. Du Val, On isolated singularities of surfaces which do not
affect the conditions of adjunction. II, III, Proc. Cambridge
Phil. Soc. 30 (1934), 483–491.
 \bibitem[HW]{HW}
F. Hidaka and K. Watanabe, Normal Gorenstein surfaces with ample
anti-canonical divisor, Tokyo J. Math. 4 (1981), no. 2, 319-330.
 \bibitem[SD]{SD}
B. Saint-Donat, Projective models of K3 surfaces, Amer. J. Math.
96 (1974), 602-639.
 \bibitem[C]{C}
F.R. Cossec, Projective models of Enriques surfaces, Math. Ann.
265 (1983), 283-334
 \bibitem[BPV]{BPV}
W. Barth, C. Peters, A. Van de Ven, Compact complex surfaces,
Heidelberg, Berlin, New York, Springer 1984.
 \bibitem[B]{B}
 E. Brieskorn, Rationale Singularit\"{a}ten Complexer Fl\"{a}chen,
 Invent. Math. 4 (1967/1968), 336-358.

 \bibitem[IF]{IF}
A.R. Iano-Fletcher, Working with weighted complete intersections,
in A.~Corti and M. Reid, editors, Explicit birational geometry of
3-folds, LMS Lecture Notes Series 281, CUP, Cambridge, 2000,
101-173.
 \bibitem[RS]{RS}
M. Reid, K. Suzuki, Cascades of projections from log del Pezzo
surfaces, in Number theory and algebraic geometry - to Peter
Swinnerton-Dyer on his 75th birthday, Cambridge University Press,
Cambridge, 2003, 227- 249.

 \bibitem[T]{T}
H.Takagi, On classification of Q-Fano 3-folds with Gorenstein
index 2 and Fano index 1/2, Preprint (1999), math.AG/9905068, 40
pages.

\end{thebibliography}
\end{document}